\documentclass[12pt]{article}
\usepackage{amssymb}
\usepackage{amsmath}
\usepackage{indentfirst,latexsym}
\usepackage{mathrsfs}
\usepackage{graphicx}
\usepackage{fancyhdr}

\date{}
\begin{document}

\title{The Eulerian distribution on self evacuated involutions}
\author{Marilena Barnabei, Flavio Bonetti, and Matteo Silimbani \thanks{
Dipartimento di Matematica - Universit\`a di Bologna}} \maketitle

\noindent {\bf Abstract.} We present an extensive study of the
Eulerian distribution on the set of self evacuated involutions,
namely, involutions corresponding to standard Young tableaux that
are fixed under the Sch$\ddot{\textrm{u}}$tzenberger map. We find
some combinatorial properties for the generating polynomial of
such distribution, together with an explicit formula for its
coefficients. Afterwards, we carry out an analogous study for the
subset of self evacuated involutions without fixed points.
\newline



\section{Introduction}

\noindent The distribution of the descent statistic (classically
known as \emph{Eulerian distribution}) on peculiar subsets of
permutations has been object of intensive studies in recent years
(see e.g. \cite{gess} and \cite{tan}) . In particular, several
authors examined the properties of the polynomial
$I_n(x)=\sum_{j=0}^{n-1}i_{n,j}x^j$, where $i_{n,j}$ denotes the
number of involutions on $[n]=\{1,2,\ldots,n\}$ with $j$ descents.
More specifically, V. Strehl \cite{stre} proved that the
coefficients of this polynomial are symmetric, and recently V.J.
Guo and J. Zeng \cite{zeng} showed that the polynomial $I_n(x)$ is
unimodal. In a previous paper \cite{bbsy} the present authors
proved that the polynomial $I_n(x)$ is not log-concave. The proof
of this property, that has been an open problem for some years,
lies upon a (not bijective) correspondence between involutions on
$[n]$ with $j$ descents and generalized involutions on length $n$
on $m$ symbols, with $m>j$. This correspondence yields an explicit
formula for the coefficients $i_{n,j}$ of the polynomial $I_n(x)$.
\newline

\noindent In this paper we study the polynomial
$S_n(x)=\sum_{j=0}^{n-1}s_{n,j}x^j$, where $s_{n,j}$ denotes the
number of \emph{self evacuated} involutions on $[n]$ with $j$
descents, namely, involutions that correspond (via the
Robinson-Schensted algorithm) to standard Young tableaux that are
fixed under the action of the Sch$\ddot{\textrm{u}}$tzenberger
map. This class of tableaux has been formerly studied by M.A.A.
van Leeuwen \cite{mvl}, who characterized the set of self
evacuated tableaux of given shape by means of domino
tilings.\newline

\noindent First of all, we exhibit an explicit formula and a
recursive rule for the total number of self evacuated involutions
on $[n]$. Following along the lines of \cite{bbsy}, we obtain some
enumerative results for the sequence $s_{n,j}$ by exploiting a map
that associates a self evacuated involution with a suitable set of
generalized involutions. In particular, we deduce an explicit
formula for the integers $s_{n,j}$, which allows to prove that the
polynomials $I_n(x)$ and $S_n(x)$ share some properties, such as
the symmetry of the coefficients and the non log-concavity.
\newline

\noindent The last section is devoted to the study of the Eulerian
distribution on self evacuated involutions without fixed points,
that is symmetric, as in the general case. Also in this case, we
find an explicit formula for the number $s^*_{n,j}$ of self
evacuated involutions on $[n]$ without fixed points and $j$ rises.

\section{Tableaux and involutions}
\label{varfy} \label{serve}

\noindent In this section, we give some definitions and general
results about tableaux, involutions and generalized
involutions.\newline

\noindent Consider the set $\mathscr{T}_n$ of standard Young
tableaux on $n$ boxes. It is well known that the
Robinson-Schensted algorithm establishes a bijection
$\rho:\mathscr{I}_n\to \mathscr{T}_n$, where $\mathscr{I}_n$ is
the set of involutions over $[n]:=\{1,2,\ldots,n\}$.\newline

\noindent We recall that the \emph{descent set} of a permutation
$\sigma$ is defined as des$(\sigma)=\{1\leq
i<n:\sigma(i)>\sigma(i+1)\}$. An analogous definition can be given
for the \emph{rise set} of a permutation, by replacing
''$\sigma(i)>\sigma(i+1)$'' by
''$\sigma(i)<\sigma(i+1)$''.\newline

\noindent Given a Ferrers diagram $\lambda$, a \emph{semistandard
tableau} of shape $\lambda$ over the alphabet $[m]$ is an array
obtained by placing into each box of the diagram $\lambda$ an
integer in $[m]$ so that the entries are strictly increasing by
rows and weakly increasing by columns.\newline

\noindent A \emph{generalized involution} is defined to be a
biword:
$$\alpha={x\choose y}=\left(\begin{array}{cccc}
x_1&x_2&\cdots&x_n\\
y_1&y_2&\cdots&y_n\\
\end{array}\right),$$
such that: \begin{itemize}\item for every $1\leq i\leq n$, there
exists an index $j$ with $x_i=y_j$ and $y_i=x_j$, \item $x_i\leq
x_{i+1}$,
\item $x_i=x_{i+1}\Longrightarrow y_i\geq
y_{i+1}$.\end{itemize} The word $x=x_1\cdots x_n$ is called the
\emph{content} of the generalized involution, and the integer $n$
is called its \emph{length}.
\newline

\noindent The Robinson-Schensted-Knuth (RSK) algorithm (see
\cite{knu}) associates bijectively a semistandard tableaux $S$
with a generalized involution inv$(S)$. \newline

\noindent We say that an integer $a$ is a \emph{repetition} of
multiplicity $r$ for the generalized involution $\alpha$ if
$$x_i=y_i=x_{i+1}=y_{i+1}=\cdots=x_{i+r-1}=y_{i+r-1}=a.$$
\newline

\noindent We define a map $\Pi$ from the set of generalized
involutions to the set of involutions as follows: if

$$\alpha=\left(\begin{array}{cccc}
x_1&x_2&\cdots&x_n\\
y_1&y_2&\cdots&y_n\\
\end{array}\right),$$
then $\Pi(\alpha)$ is the involution $\sigma$
$$\sigma=\left(\begin{array}{cccc}
1&2&\cdots&n\\
y'_1&y'_2&\cdots&y'_n\\
\end{array}\right),$$
where $y'_i=1$ if $y_i$ is the least symbol occurring in the word
$y$, $y'_j=2$ if $y_j$ is the second least symbol in $y$ and so
on. In the case $y_i=y_j$, with $i>j$, we consider $y_i$ to be
less then $y_j$. We will call the involution $\sigma=\Pi(\alpha)$
the \emph{polarization} of $\alpha$.\newline

\noindent For example, the polarization of the generalized
involution
$$\alpha=\left(\begin{array}{cccccccc}
1&1&2&3&4&4&4&6\\
4&3&2&1&6&4&1&4\\
\end{array}\right)$$
is the involution
$$\Pi(\alpha)=\left(\begin{array}{cccccccc}
1&2&3&4&5&6&7&8\\
7&4&3&2&8&6&1&5\\
\end{array}\right).$$

\noindent Note that the map $\Pi$ is not injective, since, for any
given involution $\sigma$, there are infinitely many generalized
involutions whose polarization is $\sigma$. For example, the
generalized involution

$$\beta=\left(\begin{array}{cccccccc}
1&1&1&3&4&4&5&6\\
5&3&1&1&6&4&1&4\\\end{array}\right)$$ has the same polarization as
$\alpha$ in the previous example.\newline

\noindent We will denote by Gen$_m(\sigma)$ the set of generalized
involutions, with symbols taken from $[m]$, whose polarization is
$\sigma$. Remark that two generalized involutions in
Gen$_m(\sigma)$ can not have the same content. For this reason,
the set Gen$_m(\sigma)$ corresponds bijectively with the set of
contents of its elements.

\noindent We will say that a content $x$ is \emph{compatible with}
$\sigma$ if there exists a generalized involution in some
Gen$_m(\sigma)$ whose content is $x$.

\noindent It is easy to check that a content $x=x_1\cdots x_n$ is
compatible with an involution $\sigma$ if and only if we have
$$x_i<x_{i+1}\Longleftrightarrow \sigma\textrm{ has a rise at position }i.$$

\noindent The key tool in the present paper is the interplay
between involutions and generalized involutions. For this reason,
we need to evaluate the cardinality of the set Gen$_m(\sigma)$,
for any given involution $\sigma$. It turns out that this
cardinality depends only on the number of rises of $\sigma$. In
fact, we have the following result, formerly stated in
\cite{bbsy}:

\newtheorem{yama}{Proposition}
\begin{yama}\label{nlac}
Let $\sigma\in \mathscr{I}_n$ be an involution with $t$ rises.
Then,
\begin{equation}|\textrm{Gen}_m(\sigma)|={n+m-t-1\choose n}.
\label{consalite}\end{equation}
\end{yama}

\noindent\emph{Proof} Choose an involution $\sigma\in
\mathscr{I}_n$ with $t$ rises. As we remarked above, the set
Gen$_m(\sigma)$ corresponds bijectively to the set of contents
$x=x_1\ldots x_n$ with $1\leq x_1\leq x_2\leq\cdots\leq x_{n}\leq
m$, where the inequalities are strict in correspondence of the
rises of $\sigma$. Every such content is uniquely determined by
the sequence $\delta :=\delta_0\delta_1\ldots\delta_{n}$, with
$$\delta_0=x_1-1,\quad \delta_1=x_2-x_1,\quad\ldots,\quad
\delta_n=m-x_n$$ which is a composition of the integer $m-1$ such
that its $i$-th component $\delta_i$ is at least one whenever
$\sigma$ has a rise at the $i$-th position. For this reason, we
can consider the word
$\delta'=\delta'_0\delta'_1\ldots\delta'_{n}$ defined as follows:
$$\delta'_i=\left\{\begin{array}{lc}
\delta_i-1 & \textrm{ if $\sigma$ has a rise at the i-th position}\\
\delta_i & \textrm{ otherwise}
\end{array}\right.,$$
which is a composition of the integer $m-t-1$ in $n+1$ parts. This
gives the assertion.
\begin{flushright}
$\diamond$
\end{flushright}

\section{Self evacuated standard tableaux}

\noindent We are interested in some enumerative problems
concerning Young tableaux which are fixed by the well known
Sch$\ddot{\textrm{u}}$tzenberger map (or \emph{evacuation}). First
of all we recall the definition of this map.
\newline

\noindent Given a standard Young tableau $T$ with $n$ boxes (on
$[n]$), we construct a new tableau $ev(T)$ of the same shape as
follows: we remove the symbol $1$ from the tableau $T$, leaving an
empty box. We now move into this box the smallest of the integers
contained into its two neighbor boxes. This creates a new empty
box into $T$. The process is repeated with this box according to
the same rule. It continues until there are no neighbors to slide
into the current empty box $b_1$, in which case we delete the box
$b_1$ from $T$ and we insert the symbol $n$ at the same position
in $ev(T)$. We repeat this procedure, removing from $T$ the symbol
$2$ and placing the integer $n-1$ into the box $b_2$ of $ev(T)$.
We proceed until the tableau $T$ is empty. It is well known (see
\cite{schutz}) that $ev(T)$ is a standard tableau and
$ev(ev(T))=T$.\newline

\noindent An alternative, and even simpler, description of the
Sch$\ddot{\textrm{u}}$tzenberger map can be given in terms of
involutions of the symmetric group. If $\sigma$ is the involution
associated with $T$, then the tableau $ev(T)$ corresponds to the
involution $ev(\sigma)=\psi\sigma\psi$, where $\psi$ is the
involution that maps the integer $i$ into its \emph{complement}
$n+1-i$.

\noindent This means that an involution is a fixed point under the
Sch$\ddot{\textrm{u}}$tzenberger map if and only if it is
contained in the centralizer of $\psi$. We will call such
involutions \emph{self evacuated involutions}, and the
corresponding tableaux will be called \emph{self evacuated
tableaux}.\newline

\noindent The involution point of view allows to give a simple
characterization of the fixed points of the
Sch$\ddot{\textrm{u}}$tzenberger map:\newline

\newtheorem{yama2}[yama]{Proposition}
\begin{yama2}\label{nlac}
An involution $\sigma$ is self evacuated if and only if, for every
$1\leq i\leq n$,
$$\sigma(i)+\sigma(n+1-i)=n+1.$$
\end{yama2}

\noindent \emph{Proof} The statement is a straightforward
consequence of the fact that $\sigma$ must commute with the map
$\psi$.
\begin{flushright}
$\diamond$
\end{flushright}

\noindent Recall that $\sigma$ is an involution if and only if its
disjoint cycle decomposition consists uniquely of fixed points and
transpositions. We will write $(i,j)|\,\sigma$ whenever the
transposition $(i,j)$ appears in the cycle decomposition of
$\sigma$. We will say that $(i,j)$ is a \emph{smooth
transposition} of $S_n$ if $i\neq n+1-j$. From this perspective,
Proposition \ref{nlac} can be restated as follows:

\newtheorem{eit}[yama]{Proposition}
\begin{eit}\label{lcm}
An involution $\sigma\in\mathscr{I}_n$ is self evacuated if and
only if:
\begin{equation}\sigma(i)=i\iff\sigma(n+1-i)=n+1-i,\label{fixed}\end{equation}
\begin{equation}(i,j)|\,\sigma\iff (n+1-i,n+1-j)|\,\sigma.\label{transpo}\end{equation}
\end{eit}
\begin{flushright}
$\diamond$
\end{flushright}

\noindent Note that Proposition \ref{lcm} implies that whenever a
smooth transposition divides an involution $\sigma$, this forces
four values of $\sigma$, while if a non-smooth transposition
divides $\sigma$, it forces only two values of $\sigma$.\newline

\noindent Denote by $\mathscr{S}_{n}$ the set of self evacuated
involutions on $n$ letters and by $s_n$ its cardinality.\newline

\noindent First of all, remark that $s_{2k}=s_{2k+1}$. In fact, if
$n$ is odd, Proposition \ref{nlac} implies that
$\sigma(\frac{n+1}{2})=\frac{n+1}{2}$, for every  $\sigma\in
\mathscr{S}_n$. Hence, an involution in $\mathscr{S}_{2k+1}$ is
associated to a unique involution in $\mathscr{S}_{2k}$ obtained
by deleting the central symbol.\newline

\noindent The characterization given in Proposition \ref{lcm}
allows us to give both a recurrence (Theorem \ref{orano}) and an
explicit formula (Theorem \ref{orasi}) for the integers $s_{2k}$.

\newtheorem{serec}[yama]{Theorem}
\begin{serec}
\label{orano} We have:
\begin{equation}s_{2k}=2s_{2k-2}+(2k-2)s_{2k-4}\label{ricself}\end{equation}
\end{serec}

\noindent \emph{Proof} Let $\sigma\in\mathscr{S}_{2k}$. If
$\sigma(1)=1$ or $\sigma(1)=2k$ (and hence $\sigma(2k)=2k$ or
$\sigma(2k)=1$, respectively) the restriction of $\sigma$ to the
set $\{2,\ldots,2k-1\}$ belongs to $\mathscr{S}_{2k-2}$.
Otherwise, if $\sigma(1)=j$, with $j\neq 1,2k$, we must have
$$\sigma(j)=1 \qquad \sigma(2k+1-j)=2k \qquad \sigma(2k)=2k+1-j.$$
Also in this case, the restriction of $\sigma$ to the set
$\{2,\ldots,2k-1\}\setminus\{j,2k+1-j\}$ belongs to
$\mathscr{S}_{2k-4}$. Remarking that there are $2k-2$ possible
choices for the integer $j$, we get the assertion.
\begin{flushright}
$\diamond$
\end{flushright}

\newtheorem{expself}[yama]{Theorem}
\begin{expself}
\label{orasi} The number of self evacuated involutions on $2k$
symbols is
$$s_{2k}=\sum_{h=0}^{\lfloor\frac{k}{2}\rfloor}\frac{(2k)!!}{(k-2h)!h!2^{2h}}.$$
\end{expself}

\noindent \emph{Proof} Fix an integer $h\leq
\lfloor\frac{k}{2}\rfloor$. We count the number of involutions in
$\mathscr{S}_{2k}$ with exactly $2h$ smooth transpositions. Choose
a word $w=w_1\cdots w_k$ consisting of $k$ different letters taken
from the alphabet $[2k]$ such that $w$ does not contain
simultaneously an integers $i$ and its complement $2k+1-i$. We
have $(2k)(2k-2)\cdots(2)=(2k)!!$ choices for such a word. This
word corresponds to a unique self evacuated involution $\tau$ with
$2h$ smooth transpositions defined by the following conditions:
$$\tau(w_1)=w_2,\quad\ldots\quad,\tau(w_{2h-1})=w_{2h};$$
$$\tau(w_{2h+j})=\left\{\begin{array}{lr}
w_{2h+j}&\textrm{if }w_{2h+j}\leq k\\
2k+1-w_{2h+j}&\textrm{otherwise}
\end{array}\right.,\vspace{.5cm}$$

\noindent with $0<j\leq k-2h$. It is easily checked that the
involution $\tau$ arises from $(k-2h)!h!2^{2h}$ different words
$w$. This completes the proof.
\begin{flushright}
$\diamond$
\end{flushright}

\section{Self evacuated generalized involutions}

\noindent The involution approach suggests how to extend the
Sch$\ddot{\textrm{u}}$tzenberger map to the set of semistandard
tableaux on a given alphabet $[m]$, as follows: let $S$ be a
semistandard tableau on $[m]$, with associated generalized
involution
$$\alpha=\left(\begin{array}{cccc}
x_1&x_2&\cdots&x_n\\
y_1&y_2&\cdots&y_n\\
\end{array}\right).$$ Then the
evacuated semistandard tableau $ev(S)$ is defined to be the
semistandard tableau associated with the generalized involution
$$ev(\alpha)=\left(\begin{array}{cccc}
m+1-x_n&m+1-x_{n-1}&\cdots&m+1-x_1\\
m+1-y_n&m+1-y_{n-1}&\cdots&m+1-y_1\\\end{array}\right).$$ Clearly,
the generalized involutions $\alpha$ and $ev(\alpha)$ may have
different content. More precisely, the integer $i$ occurs in the
content of $\alpha$ as many times as $m+1-i$ occurs in
$ev(\alpha)$.
\newline

\noindent For example, consider the semistandard tableau

$$S=\begin{array}{cccc}
1&2&3&4\\
2&3&4&\\
3&4&&\\
4&&&
\end{array}$$

\noindent corresponding to the generalized involution

$$\alpha=\left(\begin{array}{cccccccccc}
1&2&2&3&3&3&4&4&4&4\\
1&2&2&3&3&3&4&4&4&4\\
\end{array}\right).\vspace{.5cm}$$

\noindent The evacuated tableau is

$$ev(S)=\begin{array}{cccc}
1&2&3&4\\
1&2&3&\\
1&2&&\\
1&&&
\end{array}$$

\noindent corresponding to the generalized involution

$$ev(\alpha)=\left(\begin{array}{cccccccccc}
1&1&1&1&2&2&2&3&3&4\\
1&1&1&1&2&2&2&3&3&4\\
\end{array}\right).\vspace{.5cm}$$

\noindent From now on, extending the previous notation, we will
write $(i,j)|\alpha$ whenever the pair $(i,j)$ appears in the
generalized involution $\alpha$. Also in this case, we will say
that $(i,j)$ is a smooth transposition if $i\neq m+1-j$ and $i\neq
j$.\newline

\noindent The fixed point of the Sch$\ddot{\textrm{u}}$tzenberger
map on generalized involutions, called \emph{self evacuated
generalized involutions}, can be easily characterized as follows:

\newtheorem{format}[yama]{Proposition}
\begin{format}\label{hyu}
A generalized involution $\alpha$ is self evacuated if and only
if, whenever $(i,j)|\alpha$, we have also $(m+1-j,m+1-i)|\alpha$.
\end{format}
\begin{flushright}
$\diamond$
\end{flushright}

\noindent Remark that the Sch$\ddot{\textrm{u}}$tzenberger map
commutes with the polarization $\Pi$, namely, if $\alpha$ is a
generalized involution, we have:
$$\Pi(ev(\alpha))=ev(\Pi(\alpha)).$$

\noindent For instance, if $\alpha$ is the generalized involution
of the previous exapmle, we have:
$$\sigma_1=\Pi(\alpha)=\left(\begin{array}{cccccccccc}
1&2&3&4&5&6&7&8&9&10\\
1&3&2&6&5&4&10&9&8&7\\
\end{array}\right)$$
and
$$\sigma_2=\Pi(ev(\alpha))=\left(\begin{array}{cccccccccc}
1&2&3&4&5&6&7&8&9&10\\
4&3&2&1&7&6&5&9&8&10\\
\end{array}\right).$$
It is easily checked that $ev(\sigma_1)=\sigma_2$.\newline

\noindent Proposition \ref{hyu} yields a further characterization
of self evacuated generalized involutions, which will be useful in
the following sections.

\newtheorem{alea}[yama]{Proposition}
\begin{alea}\label{iacta}
A generalized involution $\alpha$ is self evacuated if and only if
it satisfies the following properties:
\begin{itemize}
\item the content $x=x_1\ldots x_n$ of $\alpha$ is symmetric, namely
$x_i+x_{n+1-i}=m+1,$
\item $\Pi(\alpha)$ is a self evacuated involution.
\end{itemize}
\end{alea}
\begin{flushright}
$\diamond$
\end{flushright}

\noindent We denote by $c_{n,m}$ the number of generalized
involutions of length $n$ over the alphabet $[m]$.
\newline

\noindent Setting $n=2k+1$, straightforward considerations lead to
the following properties:
\begin{itemize}
\item if $m=2h$, $c_{2k+1,m}=0$;
\item if $m=2h+1$, the central pair ${x_{k+1}\choose y_{k+1})}$ of every self evacuated generalized involution of length $n$ over the alphabet $[m]$
is necessarily the pair $(h+1,h+1)$. This implies that
$c_{2k+1,m}=c_{2k,m}$.
\end{itemize}
Hence, the values of the sequences $c_{2k+1,m}$ can derived from
the sequences $c_{2k,m}$. For this reason, we restrict to the even
case.

\newtheorem{fex}[yama]{Theorem}
\begin{fex}
The number of self evacuated generalized involutions of length
$2k$ over $[m]$ is:
\begin{equation}c_{2k,m}=\sum_{j=0}^{\lfloor\frac{k}{2}\rfloor}
{\frac{{m\choose 2}-\left\lfloor\frac{m}{2}\right\rfloor}
{2}+j-1\choose j}{m+k-2j-1\choose
k-2j}.\label{effeuno}\end{equation}
\end{fex}

\noindent \emph{Proof} Fix $h\leq\lfloor\frac{k}{2}\rfloor$. We
count the number of self evacuated generalized involutions of
length $2k$ and $m$ symbols with exactly $2h$ smooth
transpositions which, in the present case, can or can not be
different. The set $A$ of all possible smooth transposition has
cardinality
$${m\choose 2}-\left\lfloor\frac{m}{2}\right\rfloor.$$
Remark that, given a generalized involution $\alpha$ and a smooth
transposition $\tau=(i,j)$ , we have that $\tau|\,\alpha$ if and
only if $\tau'|\,\alpha$, where $\tau'=(m+1-j,m+1-i)$. It is
evident that $\tau$ can be chosen in
$$\frac{{m\choose 2}-\left\lfloor\frac{m}{2}\right\rfloor}{2}$$
ways. Such choices determine $4h$ pairs of $\alpha$. The remaining
$2k-4h$ pairs can be chosen to be either fixed points or
non-smooth transpositions. This completes the proof.
\begin{flushright}
$\diamond$
\end{flushright}

\section{The Eulerian distribution on self evacuated involutions}

\noindent In this section, we study the distribution of the
descent statistic on the set of involutions. The combinatorial
relations between involutions and generalized involutions pointed
out in the previous sections will play a crucial role for this
analysis.\newline

\noindent The distribution of the descent statistic on the set of
involutions behaves properly with respect to the action of the
Sch$\ddot{\textrm{u}}$tzenberger map. In fact:

\newtheorem{symp}[yama]{Proposition}
\begin{symp}
For every involution $\sigma$ on $[2k]$, we have:
$$|Des(\sigma)|=|Des(ev(\sigma))|.$$ Moreover, the descent sets
$Des(\sigma)$ and $Des(ev(\sigma))$ are mirror symmetric, i.e.
$\sigma$ has a descent at position $i$ if and only if $ev(\sigma)$
has a descent at position $2k-i$.
\end{symp}

\noindent \emph{Proof} Suppose that $\sigma$ has a descent at
position $i$, namely, $\sigma(i)>\sigma(i+1)$. Then,
$$ev(\sigma)(2k-i)=2k+1-\sigma(i+1)>2k+1-\sigma(i)=ev(\sigma)(2k+1-i).$$

\begin{flushright}
$\diamond$
\end{flushright}

\noindent For example, let

$$\sigma=\left(\begin{array}{cccccccc}
{\bf 1}&{\bf 2}&3&4&{\bf 5}&6&7&8\\
{\bf 3}&{\bf 2}&1&4&{\bf 6}&5&7&8\\
\end{array}\right),\vspace{.5cm}$$

\noindent where, from now on, the bold-faced numbers denote the
descent positions. Then,

$$ev(\sigma)=\left(\begin{array}{cccccccc}
1&2&{\bf 3}&4&5&{\bf 6}&{\bf 7}&8\\
1&2&{\bf 4}&3&5&{\bf 8}&{\bf 7}&6\\
\end{array}\right).\vspace{1.3cm}$$

\noindent In particular, if $\sigma$ is a self evacuated
involution, then its descent set must be mirror symmetric with
respect to the $k$-th entry.\newline

\noindent We are now interested in finding an explicit formula for
the number $s_{2k,d}$ of self evacuated involutions with $d$
rises. First of all, we have:

\newtheorem{aspe}[yama]{Proposition}
\begin{aspe}
The sequence $s_{2k,d}$ is symmetric, namely,
$$s_{2k,i}=s_{2k,2k-1-i}.$$
\end{aspe}

\noindent \emph{Proof} Given a self evacuated involution $\sigma$,
it is easily checked that the permutation $\tau=\psi\sigma$
satisfies the following properties:
\begin{itemize}
\item $\tau$ is an involution;
\item $\tau$ is self evacuated;
\item $\tau$ has a descent at position $i$ whenever $\sigma$ has a
rise at the same position.
\end{itemize}
\begin{flushright}
$\diamond$
\end{flushright}

\noindent For example, let

$$\sigma=\left(\begin{array}{cccccccc}
1&{\bf 2}&3&{\bf 4}&5&{\bf 6}&7&8\\1&{\bf 7}&5&{\bf 6}&3&{\bf
4}&2&8
\end{array}\right).$$

\noindent Then,

$$\psi\sigma=\left(\begin{array}{cccccccc}
{\bf 1}&2&{\bf 3}&4&{\bf 5}&6&{\bf 7}&8\\{\bf 8}&2&{\bf 4}&3&{\bf
6}&5&{\bf 7}&1
\end{array}\right).\vspace{1.3cm}$$

\noindent The preceding result shows that the integer $s_{2k,d}$
counts simultaneously the involutions in $\mathscr{S}_{2k}$ with
$d$ descents and those with $d$ rises.
\newline

\noindent Now we want to express the number $c_{2k,m}$ of self
evacuated generalized involutions of length $2k$ over $[m]$ in
terms of the sequence $s_{2k,d}$ by exploiting the combinatorial
relations between involutions and generalized involutions. As in
the general case (Proposition \ref{nlac}), it turns out that the
number of self evacuated generalized involutions on $m$ symbols
whose polarization is a given involution $\sigma$ depends only on
the number of rises of $\sigma$. In fact:

\newtheorem{eitdy}[yama]{Theorem}
\begin{eitdy}
\label{gilberto} We have:
\begin{equation}c_{2k,m}=\sum_{j=0}^{m-1}{k+\left\lfloor\frac{j}{2}\right\rfloor\choose
\left\lfloor\frac{j}{2}\right\rfloor}s_{2k,m-1-j}.\label{effedue}\end{equation}
\end{eitdy}

\noindent \emph{Proof} Let $\sigma\in \mathscr{I}_{2k}$ a self
evacuated involution with $t$ rises. As remarked in proposition
\ref{nlac}, $\sigma$ corresponds to
$${2k+m-1-t\choose m-1-t}$$
generalized involutions with $m$ symbols, but only
$${k+\left\lfloor\frac{m-1-t}{2}\right\rfloor\choose
\left\lfloor\frac{m-1-t}{2}\right\rfloor}$$ of these are self
evacuated. In fact, by Proposition \ref{iacta}, a generalized
involution with $m$ symbols in the set Gen$_m(\sigma)$ is self
evacuated if only if the corresponding composition $\delta'$ of
the integer $m-1-t$ into $2k+1$ satisfies the condition
$\delta'_{k-i}=\delta'_{k+i}$. By setting $j=m-1-t$, we get the
assertion.
\begin{flushright}
$\diamond$
\end{flushright}

\noindent We now exploit the described combinatorial relation
between generalized involutions and involutions to determine an
explicit formula for $s_{2k,d}$.

\newtheorem{finalf}[yama]{Theorem}
\begin{finalf}
The number of self evacuated involutions of length $2k$ with $d$
rises is:
\begin{equation}s_{2k,d}=\sum_{j=1}^{d+1}(-1)^{\left\lfloor\frac{d-j}{2}+1\right\rfloor}{k\choose\left\lfloor\frac{d+1-j}{2}\right\rfloor}
\sum_{i=0}^{\lfloor\frac{k}{2}\rfloor} {\frac{{j\choose
2}-\left\lfloor\frac{j}{2}\right\rfloor} {2}+i-1\choose
i}{j+k-2i-1\choose k-2i}.\label{effefin}\end{equation}
\end{finalf}

\noindent \emph{Proof} Formula (\ref{effedue}) yields, by
inversion:
\begin{equation}s_{2k,d}=\sum_{j=1}^{d+1}(-1)^{\left\lfloor\frac{d-j}{2}+1\right\rfloor}{k\choose\left\lfloor\frac{d+1-j}{2}\right\rfloor}
c_{2k,j}.\label{effetre}\end{equation} Then, combining Formulae
(\ref{effeuno}) and (\ref{effetre}), we derive (\ref{effefin}).
\begin{flushright}
$\diamond$
\end{flushright}

\noindent Moreover, this explicit formula allows to check that the
polynomials $S_{2k}(x)=\sum_{j=0}^{2k-1}s_{2k,j}x^j$ are not, in
general, log-concave, since we have, for example:
$$s_{100,0}\cdot s_{100,2}=11950>2500=s_{100,1}^2.$$\newline

\noindent The first values of $s_{2k,d}$ are shown in the
following table:

$$\begin{array}{l|llllllllll}
n/d&0&1&2&3&4&5&6&7&8&9\\\hline
0&1&&&&&&&&&\\
1&1&&&&&&&&&\\
2&1&1&&&&&&&&\\
3&1&0&1&&&&&&&\\
4&1&2&2&1&&&&&&\\
5&1&0&4&0&1&&&&&\\
6&1&3&6&6&3&1&&&&\\
7&1&0&9&0&9&0&1&&&\\
8&1&4&13&20&20&13&4&1&&\\
9&1&0&17&0&40&0&17&0&1&\\
10&1&5&23&49&78&78&49&23&5&1\\
\end{array}$$

\noindent These first values seem to suggest that the polynomials
$S_{2k}(x)$ are unimodal for every $k\in\mathbb{N}$. It would be
interesting to find a combinatorial proof of this property.

\section{Self evacuated involutions without fixed points}

\noindent In this section, we extend the study of the Eulerian
distribution to the set of self evacuated involutions on $[n]$
 without fixed points. Obviously, such involutions exist
only if $n$ is even.

\noindent Denote by $\mathscr{S}^*_{2k}$ the set of self evacuated
involutions on $2k$ objects without fixed points and by $s^*_{2k}$
the cardinality of $\mathscr{S}^*_{2k}$. Then:

\newtheorem{krkr}[yama]{Theorem}
\begin{krkr}
We have:
\begin{equation}s^*_{2k}=\sum_{h=0}^{\left\lfloor\frac{k}{2}\right\rfloor}\frac{k!}{(k-2h)!h!},\label{nofix}\end{equation}
and
\begin{equation}s^*_{2k}=s^*_{2k-2}+(2n-2)s^*_{2k-4}.\label{nofixdue}\end{equation}
\end{krkr}

\noindent \emph{Proof} Following along the lines of the proof of
Theorem \ref{orasi}, we count the number of self evacuated
involutions without fixed points with exactly $2h$ smooth
transpositions, $2h\leq k$. Choose a word $w=w_1\cdots w_k$
consisting of $k$ different letters taken from the alphabet
${1,\ldots,2k}$ such that $w$ does not contain simultaneously the
integers $i$ and $2k+1-i$. We have $(2k)!!$ choices for such a
word. This word corresponds to a unique self evacuated involution
$\tau$ without fixed points with $2h$ smooth transpositions
defined by the following conditions:
$$\tau(w_1)=w_2,\quad\ldots,\quad\tau(w_{2h-1})=w_{2h},$$
$$\tau(w_{2h+j})=
2k+1-w_{2h+j},\quad\textrm{ for }0<j\leq k-2h.$$

\noindent It is easily checked that the involution $\tau$ arises
from $(k-2h)!h!2^k$ different words $w$. Hence:
$$s^*_{2k}=\sum_{h=0}^{\left\lfloor\frac{k}{2}\right\rfloor}\frac{(2k)!!}{(k-2h)!h!2^k},$$
which is equivalent to (\ref{nofix}).\newline

\noindent Let now $\sigma\in\mathscr{S}^*_{2k}$. If
$\sigma(1)=2k$, and hence $\sigma(2k)=1$, the restriction of
$\sigma$ to the set $\{2,\ldots,2k-1\}$ is a self evacuated
involution on $2k-2$ symbols without fixed points. If
$\sigma(1)=j$, with $j<2k$, the symbol $1$ is involved in a smooth
transposition, hence we must have $\sigma(j)=1$,
$\sigma(2k+1-j)=2k$ and $\sigma(2k)=2k+1-j$. Then, the restriction
of $\sigma$ to the set $\{2,\ldots,2k-1\}\setminus\{j,2k+1-j\}$ is
a self evacuated involution on $2k-4$ symbols without fixed
points. Remarking that there are $2k-2$ possible choices for the
integer $j$, we get (\ref{nofixdue}).
\begin{flushright}
$\diamond$
\end{flushright}

\noindent Denote by $s^*_{2k,d}$ the number of involutions in
$\mathscr{S}^*_{2k}$ with $d$ rises. Then:

\newtheorem{sy2}[yama]{Proposition}
\begin{sy2}
The sequence $s^*_{2k,d}$ is symmetric, namely,
$$s^*_{2k,d}=s^*_{2k,2k-d}.$$
\end{sy2}

\noindent \emph{Proof} Denote by $\mathscr{I}^*_{2k}$ the set of
involutions on $2k$ objects without fixed points. In \cite{stre},
V. Strehl proved the symmetry of the Eulerian distribution on
$\mathscr{I}^*_{2k}$ by means of a bijection:
$$\theta:\mathscr{I}^*_{2k}\to \mathscr{I}^*_{2k},$$
which maps an involutions $\sigma$ with $j$ rises to an
involutions $\theta(\sigma)$ with $2k-j$ rises. It is easily
checked that the restriction of $\theta$ to the set
$\mathscr{S}^*_{2k}$ is a bijections of $\mathscr{S}^*_{2k}$ into
itself. This gives the assertion.
\begin{flushright}
$\diamond$
\end{flushright}

\noindent Once more, in order to find
 an explicit formula for
the integers $s^*_{2k,d}$, we need to establish a connection
between self evacuated involutions without fixed points and a
suitable set of generalized involutions. Remarked that this set
contains only self evacuated generalized involutions with
repetitions of even multiplicity. Denote by $c^*_{2k,m}$ the
number of such involutions of length $2k$ on the alphabet $[m]$.
Then:

\newtheorem{xef}[yama]{Theorem}
\begin{xef}
We have:
\begin{equation}c^*_{2k,m}=\sum_{j=0}^{\lfloor\frac{k}{2}\rfloor}
{\frac{{m\choose 2}+\left\lfloor\frac{m}{2}\right\rfloor}
{2}+j-1\choose j}{\left\lceil\frac{m}{2}\right\rceil+k-2j-1\choose
k-2j}.\label{effequalche}\end{equation}
\end{xef}

\noindent \emph{Proof} Remark that, given a generalized involution
$\sigma$ and a smooth transposition $\tau=(i\ \,j)$ , we have that
$(i\ \,j)|\,\sigma$ if and only if $\tau'=(m+1-j\
\,m+1-i)|\,\sigma$. Similarly, every \emph{non central fixed
point}, namely, an occurrence of a pair $(i\ \,i)$ in $\sigma$,
with $i\neq \frac{m+1}{2}$, implies a second occurrence of the
same pair.

\noindent Fix now $j\leq\lfloor\frac{k}{2}\rfloor$. We count the
number of generalized involutions of length $2k$ on the alphabet
$[m]$ containing only repetitions of even multiplicity, such that
exactly $2j$ of its pairs are either non central fixed points or
smooth transpositions. We can choose a non central fixed point in
$\left\lfloor\frac{m}{2}\right\rfloor$ ways and a smooth
transposition in $\frac{{m\choose
2}+\left\lfloor\frac{m}{2}\right\rfloor}{2}$ ways. The remaining
pairs must be chosen to be either cental fixed points or a non
smooth transpositions. This completes the proof.
\begin{flushright}
$\diamond$
\end{flushright}

\noindent Repeating the same argumentations as in the proof of
Theorem \ref{gilberto}, we obtain the following result:

\newtheorem{edy}[yama]{Theorem}
\begin{edy}
We have:
\begin{equation}c^*_{2k,m}=\sum_{j=0}^{m-1}{k+\left\lfloor\frac{j}{2}\right\rfloor\choose
\left\lfloor\frac{j}{2}\right\rfloor}s^*_{2k,m-1-j}.\label{eedue}\end{equation}
Hence:
\begin{equation}s^*_{2k,d}=\sum_{j=1}^{d+1}(-1)^{\left\lfloor\frac{d-j}{2}+1\right\rfloor}{k\choose\left\lfloor\frac{d+1-j}{2}\right\rfloor}
\sum_{i=0}^{\lfloor\frac{k}{2}\rfloor} {\frac{{j\choose
2}+\left\lfloor\frac{j}{2}\right\rfloor} {2}+i-1\choose
i}{\left\lceil\frac{j}{2}\right\rceil+k-2i-1\choose
k-2i}.\label{eefin}\end{equation}
\end{edy}
\begin{flushright}
$\diamond$
\end{flushright}

\noindent The present table contains the first values of the
sequences $s^*_{2k,d}$:

$$\begin{array}{l|llllllllll}
n/d&0&1&2&3&4&5&6&7&8&9\\\hline
0&1&&&&&&&&&\\
2&1&&&&&&&&&\\
4&1&1&1&&&&&&&\\
6&1&1&3&1&1&&&&&\\
8&1&2&7&5&7&2&1&&&\\
10&1&2&12&12&27&12&12&2&1&\\
\end{array}$$

\noindent This table shows that the polynomial $S^*_{2k}(x)$ is
not in general unimodal, and hence not log-concave.


\begin{thebibliography}{99}

\bibitem{bbsy} M.Barnabei, F.Bonetti, M.Silimbani, The descent statistic on involutions is not
log-concave, to appear

\bibitem{gess} I.M.Gessel, C.Reutenauer, Counting permutations
with a given cycle structure and descent set, \emph{J. Combin.
Theory Ser. A} {\bf 13} (1972), 135-139.

\bibitem{zeng} V.J.Guo, J.Zeng, The Eulerian distribution on involutions is indeed unimodal,
 \emph{J. Combin. Theory Ser. A} {\bf 113} (2006), no. {\bf 6}, 1061--1071.

\bibitem{knu}
D.E.Knuth, Permutations, Matrices and Generalized Young Tableaux,
\emph{Pacific J. Math.} \textbf{34} (1970), 709-727.

\bibitem{schutz} M.P.Sch$\ddot{\textrm{u}}$tzenberger, Quelques Remarques sur une Construction de Schensted, \emph{Math. Scand.}
\textbf{12} (1963), 117-128.

\bibitem{stabook} R.P.Stanley, Enumerative Combinatorics, Vol. II, \emph{Cambridge
Studies in Advanced Mathematics}, {\bf 62}. Cambridge University
Press, Cambridge (1999).

\bibitem{stre} V.Strehl, Symmetric Eulerian distributions for
involutions, \emph{S\'eminaire Lotharingien Combinatoire} {\bf 1},
Strasbourg 1980, Publications del l'I.R.M.A. 140/S-02, Strasbourg
1981.

\bibitem{tan} S.Tanimoto, A study of Eulerian numbers for
permutations in the alternating group, \emph{Integers}, {\bf 6},
A31 (2006) (electronic)

\bibitem{mvl} M.A.A. van Leuween, The Robinson-Schensted and Sch$\ddot{\textrm{u}}$tzenberger algorithms, an elementary
approach, \emph{Electron. J. Combin.}, {\bf 3}, No. 2 (1996),
391-422.

\end{thebibliography}
\end{document}